\documentclass[a4paper,10pt]{article}

\usepackage[a4paper,left=2cm,right=1.25cm, top=2cm]{geometry}
\usepackage[english]{babel}
\usepackage{fancyhdr}
\usepackage[dvips]{graphicx}
\usepackage[latin1]{inputenc}
\usepackage{picins}
\usepackage{latexsym}
\usepackage[]{subfigure}
\usepackage{booktabs}
\usepackage[T1]{fontenc}
\usepackage{times}
\usepackage{multicol}

\usepackage{url}

\usepackage{amsmath,amssymb,amsthm}

\newtheorem{theorem}{Theorem}
\newtheorem{proposition}[theorem]{Proposition}
\newtheorem{lemma}[theorem]{Lemma}
\newtheorem{corollary}[theorem]{Corollary}
\newtheorem{example}{Example}
\newtheorem{definition}[theorem]{Definition}
\newtheorem{remark}{Remark}

\begin{document}

\makeatletter
\newenvironment{tablehere}
  {\def\@captype{table}}
  {}

\newenvironment{figurehere}
  {\def\@captype{figure}}
  {}
\makeatother

\renewcommand{\figurename}{\textbf{Fig.}}


\noindent {\LARGE\textbf{The Variational Calculus on Time Scales}}\\

\vspace{0.5cm}

\noindent \small
Delfim F. M. Torres\textsuperscript{1,}\footnote{Corresponding author: delfim@ua.pt}\\
\noindent {\textsuperscript{1}Department of Mathematics,
University of Aveiro, 3810-193 Aveiro, Portugal}

\begin{quotation}
\noindent {\textbf{Abstract.}} The discrete, the quantum,
and the continuous calculus of variations,
have been recently unified and extended
by using the theory of time scales.
Such unification and extension is,
however, not unique, and two approaches
are followed in the literature: one dealing with minimization
of delta integrals; the other dealing with minimization
of nabla integrals. Here we review a more general approach to the
calculus of variations on time scales that allows to obtain
both delta and nabla results as particular cases.

\vspace{0.5cm}

\end{quotation}

\vspace{0.5cm}

\normalsize

\begin{multicols}{2}


\section{Introduction}

The calculus on time scales is a recent field introduced by B.~Aulbach and S.~Hilger
in order to unify the theories of difference and differential equations \cite{A:H,SH}.
It has found applications in several contexts that require simultaneous modeling
of discrete and continuous data, and is nowadays under strong current research
in areas as diverse as control of population, quantum calculus, economics,
communication networks and robotic control (see \cite{A:B:O:P:02,SSW}
and references therein). The area of the calculus of variations on time scales,
which we are concerned here, is in its beginning and is a
fertile area of research \cite{F:T:07}. As particular cases,
one gets the classical calculus of variations \cite{GelfandFomin},
the discrete-time calculus of variations \cite{KP},
and the $q$-calculus of variations \cite{q:B:04}.

The calculus of variations on time scales was introduced
in 2004 by M.~Bohner using the delta derivative and integral
\cite{B:04}, and has been since then further developed
by several different authors in several different directions
\cite{A:M:02,B:T:08,F:T:08,F:T:09,AM:T,Basia:post_doc_Aveiro:1}.
In all these works, the integral functional to be extremized
has the form
\begin{equation}
\label{eq:Pd}
\mathcal{J}_\Delta(y)
= \int_a^b L\left(t,y^\sigma(t),y^\Delta(t)\right) \Delta t \, .
\end{equation}
Motivated by applications in economics \cite{A:B:L:06,A:U:08},
a different formulation for the problems of the calculus
of variations on time scales has been considered, which
involve a functional with a nabla derivative and
a nabla integral \cite{A:T,Atici:comparison,NM:T}:
\begin{equation}
\label{eq:Pn}
\mathcal{J}_\nabla(y)
= \int_a^b L\left(t,y^\rho(t),y^\nabla(t)\right) \nabla t \, .
\end{equation}
Formulations \eqref{eq:Pd} and \eqref{eq:Pn} are good
in the sense that results obtained \emph{via}
delta and nabla approaches are similar among them
and similar to the classical results of the calculus
of variations. An example of this is given
by the time scale versions of the Euler--Lagrange equations:
if $y \in C_{\textrm{rd}}^2$ is an extremizer of \eqref{eq:Pd},
then $y$ satisfies the delta-differential equation
\begin{equation}
\label{eq:el:d} \frac{\Delta}{\Delta t}
\partial_{3}L\left(t,y^\sigma(t),{y}^\Delta(t)\right)
= \partial_{2}L\left(t,y^\sigma(t),{y}^\Delta(t)\right)
\end{equation}
for all $t \in [a,b]^{\kappa^2}$ \cite{B:04};
if $y \in C_{\textrm{ld}}^2$ is an extremizer of \eqref{eq:Pn},
then $y$ satisfies the nabla-differential equation
\begin{equation}
\label{eq:el:n} \frac{\nabla}{\nabla t}
\partial_{3}L\left(t,y^\rho(t),{y}^\nabla(t)\right)
= \partial_{2}L\left(t,y^\rho(t),{y}^\nabla(t)\right)
\end{equation}
for all $t \in [a,b]_{\kappa^2}$ \cite{NM:T},
where we use $\partial_{i}L$ to denote the standard
partial derivative of $L(\cdot,\cdot,\cdot)$
with respect to its $i$th variable, $i = 1,2,3$.
In the classical context $\mathbb{T} = \mathbb{R}$
one has
\begin{equation}
\label{eq:Pc}
\mathcal{J}_\Delta(y)
= \mathcal{J}_\nabla(y)
= \int_a^b L\left(t,y(t),y'(t)\right) dt
\end{equation}
and both \eqref{eq:el:d} and \eqref{eq:el:n} coincide
with the standard Euler--Lagrange equation:
if $y \in C^2$ is an extremizer of
the integral functional \eqref{eq:Pc}, then
\begin{equation*}
\frac{d}{d t}
\partial_{3}L\left(t,y(t),y'(t)\right)
=\partial_{2}L\left(t,y(t),y'(t)\right)
\end{equation*}
for all $t \in [a,b]$.
However, the problems of extremizing \eqref{eq:Pd}
and \eqref{eq:Pn} are intrinsically different, in the sense
that is not possible to obtain the nabla results as corollaries
of the delta ones and \emph{vice versa}. Indeed,
if admissible functions $y$ are of class $C^2$
then (\textrm{cf.} \cite{G:G:S:05})
\begin{equation*}
\begin{split}
\mathcal{J}_\Delta(y)
&= \int_a^b L\left(t,y^\sigma(t),y^\Delta(t)\right) \Delta t\\
&= \int_a^b L\left(\rho(t),(y^\sigma)^\rho(t),y^\nabla(t)\right) \nabla t
\end{split}
\end{equation*}
while
\begin{equation*}
\begin{split}
\mathcal{J}_\nabla(y)
&= \int_a^b L\left(t,y^\rho(t),y^\nabla(t)\right) \nabla t\\
&= \int_a^b L\left(\sigma(t),(y^\rho)^\sigma(t),y^\Delta(t)\right) \Delta t
\end{split}
\end{equation*}
and one easily see that functionals
\eqref{eq:Pd} and \eqref{eq:Pn}
have a different nature and are not
compatible with each other.

In this paper we consider the more general delta-nabla formulation
of the calculus of variations introduced in \cite{delfim:Bedlewo:2009}
and further developed in \cite{china-Xuzhou,MyID:180,china-Xiamen,MOMA09},
that includes, as trivial examples, the problems with functionals
$\mathcal{J}_\Delta(y)$ and $\mathcal{J}_\nabla(y)$ that have been
previously studied in the literature of time scales.
For a different approach, based on the concept of duality \cite{Caputo},
we refer the reader to \cite{comBasia:duality,MyID:174,MyID:177}.
One of the main results provide an Euler--Lagrange
necessary optimality type condition
for the more general integral functional $\mathcal{J}$
we are promoting here (Section~\ref{subsec:mr:el}).
As straight corollaries we obtain
both delta- and nabla-differential equations \eqref{eq:el:d}
and \eqref{eq:el:n}. Natural boundary conditions and necessary
optimality conditions for delta-nabla isoperimetric problems
of the calculus of variations are obtained as well
(Sections \ref{subsec:mr:nbc} and \ref{subsec:mr:iso}, respectively).
Simple illustrative examples on the application of the results
are given in detail.


\section{The goal}

Let $\mathbb{T}$ be a given time scale with $a, b \in \mathbb{T}$, $a < b$,
and $\left(\mathbb{T} \setminus \{a,b\}\right)\cap [a,b] \ne \emptyset$;
$L_{\Delta}(\cdot,\cdot,\cdot)$ and $L_{\nabla}(\cdot,\cdot,\cdot)$ be two given smooth
functions from $\mathbb{T} \times \mathbb{R}^2$ to $\mathbb{R}$.
The results here discussed are trivially generalized for
admissible functions $y : \mathbb{T}\rightarrow\mathbb{R}^n$
but for simplicity of presentation
we restrict ourselves to the scalar case $n=1$.
We consider the delta-nabla integral functional
\begin{equation}
\label{eq:P}
\begin{split}
\mathcal{J}(y)
&= \left(\int_a^b L_{\Delta}[y](t) \Delta t\right)
\cdot
\left(\int_a^b L_{\nabla}\{y\}(t) \nabla t\right)\\
&= \int_a^b \int_a^b
L_{\Delta}[y](t) \cdot L_{\nabla}\{y\}(t) \Delta t \nabla \tau
\end{split}
\end{equation}
where, for simplicity of notation, we use the operators
$[\cdot]$ and $\{\cdot\}$ defined by
\begin{equation*}
\begin{split}
[y](t) &= \left(t,y^\sigma(t),y^\Delta(t)\right) \, , \\
\{y\}(t) &= \left(t,y^\rho(t),y^\nabla(t)\right) \, .
\end{split}
\end{equation*}
Note that with the new operators
we write \eqref{eq:Pd} and \eqref{eq:Pn} as
\begin{equation*}
\begin{split}
\mathcal{J}_\Delta(y) &= \int_a^b L_\Delta[y](t) \Delta t \, , \\
\mathcal{J}_\nabla(y) &= \int_a^b L_\nabla\{y\}(t) \nabla t \, .
\end{split}
\end{equation*}

\begin{remark}
\label{obs}
In the particular case $L_\nabla \equiv \frac{1}{b-a}$
functional \eqref{eq:P} reduces to \eqref{eq:Pd}
(\textrm{i.e.}, $\mathcal{J}(y) = \mathcal{J}_\Delta(y)$);
in the particular case $L_\Delta \equiv \frac{1}{b-a}$
functional \eqref{eq:P} reduces to \eqref{eq:Pn}
(\textrm{i.e.}, $\mathcal{J}(y) = \mathcal{J}_\nabla(y)$).
\end{remark}

\medskip

Our main goal is to answer the following question:
\emph{What is the Euler--Lagrange equation for
$\mathcal{J}(y)$ defined by \eqref{eq:P}?}

\medskip


\fancyhf{} \pagestyle{fancy}
\chead{D.F.M. Torres: The variational calculus on time scales.}


\section{The time scales calculus}

The aim of the calculus on time scales is to unify continuous and discrete analysis
into a general theory. The motivation for such general theory is
rooted in the fact that many results concerning differential equations
carry over quite easily to corresponding results for difference equations,
while other results seem to be very different from continuous counterparts.
The unification and extension given by the theory of time scales
helps to explain such similarities and discrepancies.

In this section we introduce the basic definitions and results
that will be needed in the sequel. For a more general presentation
of the theory of time scales and detailed proofs,
we refer the reader to the books \cite{B:P:01,B:P:03,Lak:book}.

As usual,  $\mathbb{R}$, $\mathbb{Z}$,
and $\mathbb{N}$ denote, respectively, the set of real,
integer, and natural numbers.


\subsection{The delta calculus}

A {\it time scale} $\mathbb{T}$ is an arbitrary nonempty closed subset
of $\mathbb{R}$.  Besides standard cases of $\mathbb{R}$ (continuous time)
and $\mathbb{Z}$ (discrete time), many different models of time may be used,
\textrm{e.g.}, the $h$-numbers ($\mathbb{T}$ = $h\mathbb{Z}:=\{h z \ | \  z \in
\mathbb{Z}\}$, where $h>0$ is a fixed real number) and the
$q$-numbers ($\mathbb{T}$ = $q^{\mathbb{N}_0}:=\{q^k \ | \  k \in \mathbb{N}_0\}$,
where $q>1$ is a fixed real number). We assume that a time scale
$\mathbb{T}$ has the topology that it inherits from the real numbers
with the standard topology.
For each time scale $\mathbb{T}$
the following operators are used:

\begin{itemize}
\item the {\it forward jump operator} $\sigma:\mathbb{T} \rightarrow \mathbb{T}$,
defined by $\sigma(t):=\inf\{s \in \mathbb{T}:s>t\}$ for $t\neq\sup \mathbb{T}$
and $\sigma(\sup\mathbb{T})=\sup\mathbb{T}$ if $\sup\mathbb{T}<+\infty$;

\item the {\it backward jump operator} $\rho:\mathbb{T} \rightarrow \mathbb{T}$,
defined by $\rho(t):=\sup\{s \in \mathbb{T}:s<t\}$ for $t\neq\inf \mathbb{T}$
and $\rho(\inf\mathbb{T})=\inf\mathbb{T}$ if $\inf\mathbb{T}>-\infty$;

\item the {\it forward graininess function} $\mu:\mathbb{T} \rightarrow [0,\infty)$,
defined by $\mu(t):=\sigma(t)-t$;

\item the \emph{backward graininess function}
$\nu:\mathbb{T}\rightarrow[0,\infty)$, defined by
$\nu(t)=t - \rho(t)$.

\end{itemize}

\begin{example}
If $\mathbb{T}=\mathbb{R}$, then for any $t \in \mathbb{T}$
$\sigma(t)=\rho(t)=t$ and $\mu(t) = \nu(t) = 0$.
If $\mathbb{T}=h\mathbb{Z}$, $h > 0$, then for every $t \in \mathbb{T}$
$\sigma(t)=t+h$, $\rho(t)=t-h$, and $\mu(t) = \nu(t) = h$.
On the other hand, if $\mathbb{T} = q^{\mathbb{N}_0}$, $q>1$,
then we have $\sigma(t) = q t$, $\rho(t) = q^{-1} t$,
$\mu(t)= (q-1) t$, and $\nu(t)= (1-q^{-1}) t$.
\end{example}

A point $t\in\mathbb{T}$ is called \emph{right-dense},
\emph{right-scattered}, \emph{left-dense} or \emph{left-scattered}
if $\sigma(t)=t$, $\sigma(t)>t$, $\rho(t)=t$,
and $\rho(t)<t$, respectively. We say that $t$ is \emph{isolated}
if $\rho(t)<t<\sigma(t)$, that $t$ is \emph{dense} if $\rho(t)=t=\sigma(t)$.

If $\sup \mathbb{T}$ is finite and left-scattered, we define
$\mathbb{T}^\kappa := \mathbb{T}\setminus \{\sup\mathbb{T}\}$,
otherwise $\mathbb{T}^\kappa :=\mathbb{T}$.

\begin{definition}
\label{def:de:dif} We say that a function
$f:\mathbb{T}\rightarrow\mathbb{R}$ is \emph{delta differentiable}
at $t\in\mathbb{T}^{\kappa}$ if there exists a number
$f^{\Delta}(t)$ such that for all $\varepsilon>0$ there is a
neighborhood $U$ of $t$ such that
$$
|f(\sigma(t))-f(s)-f^{\Delta}(t)(\sigma(t)-s)|
\leq\varepsilon|\sigma(t)-s|
$$
for all $s\in U$. We call
$f^{\Delta}(t)$ the \emph{delta derivative} of $f$ at $t$ and $f$ is
said \emph{delta differentiable} on $\mathbb{T}^{\kappa}$ provided
$f^{\Delta}(t)$ exists for all $t\in\mathbb{T}^{\kappa}$.
\end{definition}

\begin{remark}
If $t \in \mathbb{T} \setminus \mathbb{T}^\kappa$, then
$f^{\Delta}(t)$ is not uniquely defined, since for such a point $t$,
small neighborhoods $U$ of $t$ consist only of $t$ and, besides, we
have $\sigma(t) = t$. For this reason, maximal left-scattered points
are omitted in Definition~\ref{def:de:dif}.
\end{remark}

Note that in right-dense points
$$
f^{\Delta} (t)=\lim_{s\rightarrow t}\frac{f(t)-f(s)}{t-s},
$$
provided this limit exists, and in
right-scattered points
$$
f^{\Delta}(t)=\frac{f(\sigma(t))-f(t)}{\mu(t)},
$$
provided $f$ is continuous at $t$.

\begin{example}
\label{ex:der:QC}
If $\mathbb{T}=\mathbb{R}$, then $f:\mathbb{R} \rightarrow \mathbb{R}$ is
delta differentiable at $t \in \mathbb{R}$ if and only if $f$ is differentiable in
the ordinary sense at $t$. Then, $f^{\Delta}(t)=f'(t)$.
If $\mathbb{T} = h\mathbb{Z}$, $h > 0$, then
$f:\mathbb{Z} \rightarrow \mathbb{R}$ is always delta differentiable
at every $t \in \mathbb{Z}$ with $f^{\Delta}(t) = \frac{f(t+h) - f(t)}{h}$.
If $\mathbb{T} = q^{\mathbb{N}_0}$, $q>1$,
then $f^{\Delta} (t)=\frac{f(q t)-f(t)}{(q-1) t}$,
\textrm{i.e.}, we get the usual $q$-derivative
of quantum calculus \cite{QC}. For a more general
quantum calculus see \cite{withMiguel01,MyID:187}.
\end{example}

Delta derivatives of higher order are defined in
the standard way. Let $r\in\mathbb{N}$,
$\mathbb{T}^{\kappa^{0}} := \mathbb{T}$, and
$\mathbb{T}^{\kappa^i}:=\left(\mathbb{T}^{\kappa^{i-1}}\right)^\kappa$,
$i = 1, \ldots, r$. For convenience we also put $f^{\Delta^0} = f$
and $f^{\Delta^1} = f^\Delta$. We define the $r^{th}$-\emph{delta derivative}
$f^{\Delta^r}$ of $f$ provided $f^{\Delta^{r-1}}$ is delta differentiable
on $\mathbb{T}^{\kappa^r}$ with derivative $f^{\Delta^r} =
\left(f^{\Delta^{r-1}}\right)^\Delta: \mathbb{T}^{\kappa^r} \rightarrow \mathbb{R}$.

We shall often denote $f^\Delta$ by $\frac{\Delta}{\Delta t} f$
if $f$ is a composition of other functions of $t$.
For $f:\mathbb{T} \rightarrow X$, where $X$ is an arbitrary set,
we define $f^\sigma:=f\circ\sigma$.

For delta differentiable $f$ and $g$, the next formulas hold:

\begin{align*}
f^\sigma(t)&=f(t)+\mu(t)f^\Delta(t) \, ,\\
(fg)^\Delta(t)&=f^\Delta(t)g^\sigma(t)+f(t)g^\Delta(t)\\
&=f^\Delta(t)g(t)+f^\sigma(t)g^\Delta(t).
\end{align*}

Let $a,b \in \mathbb{T}$, $a<b$.
We define the interval $[a,b]$ in $\mathbb{T}$ by
$$[a,b]:=\{ t \in \mathbb{T}: a\leq t\leq b\}.$$
Open intervals, half-open intervals and unbounded intervals
in $\mathbb{T}$ are defined accordingly.

In order to describe a class of functions that possess a delta
antiderivative, the following definition is introduced:

\begin{definition}
A function $f:\mathbb{T}\rightarrow\mathbb{R}$ is called
\emph{rd-continuous} if it is continuous at right-dense points
in $\mathbb{T}$ and its left-sided limits exist (finite)
at left-dense points in $\mathbb{T}$. We denote the
set of all rd-continuous functions by
$C^0_{\textrm{rd}} = C_{\textrm{rd}}
= C_{\textrm{rd}}(\mathbb{T})
= C_{\textrm{rd}}(\mathbb{T}; \mathbb{R})$.
The set of functions $f : \mathbb{T} \rightarrow \mathbb{R}$
that are delta differentiable and whose delta derivative
is rd-continuous is denoted by $C_{\textrm{rd}}^1
= C_{\textrm{rd}}^1(\mathbb{T})
= C_{\textrm{rd}}^1(\mathbb{T}; \mathbb{R})$.
In general, we say that $f \in C_{\textrm{rd}}^r$, $r \in \mathbb{N}$,
if $f^\Delta \in C_{\textrm{rd}}^{r-1}$.
\end{definition}

\begin{definition}
A function $F:\mathbb{T}\rightarrow\mathbb{R}$ is called a
\emph{delta antiderivative} of
$f:\mathbb{T}\rightarrow\mathbb{R}$ provided
$$
F^{\Delta}(t)=f(t), \qquad  \forall t \in \mathbb{T}^\kappa.
$$
In this case we define the \emph{delta integral} of $f$ from $a$
to $b$ ($a,b \in \mathbb{T}$) by
\begin{equation*}
\int_{a}^{b}f(t)\Delta t:=F(b)-F(a) \, .
\end{equation*}
\end{definition}

\begin{theorem}[Theorem~1.74 of \cite{B:P:01}]
Every rd-continuous function has a delta
antiderivative. In particular, if $a \in \mathbb{T}$,
then the function $F$ defined by
$$
F(t)= \int_{a}^{t}f(\tau)\Delta\tau, \quad t \in \mathbb{T} \, ,
$$
is a delta antiderivative of $f$.
\end{theorem}

\begin{example}
Let $a, b \in \mathbb{T}$ with $a < b$
and $f \in C_{\textrm{rd}}$. (i) If $\mathbb{T} =
\mathbb{R}$, then $\int_{a}^{b}f(t)\Delta t = \int_{a}^{b}f(t) dt$,
where the integral on the right-hand side is the classical Riemann
integral. (ii) If $\mathbb{T}=h\mathbb{Z}$, $h>0$, then
\[
\int\limits_{a}^{b} f(t) \Delta
t=\sum\limits_{k=\frac{a}{h}}^{\frac{b}{h}-1} f(kh) h\, .
\]
If $\mathbb{T} = q^{\mathbb{N}_0}$, $q>1$, then
$$
\int_{a}^{b}f(t)\Delta t = (1 - q) \sum_{t \in [a,b)} f(t) t
$$
(\textrm{cf.} the usual $q$-integral of quantum calculus \cite{QC}).
\end{example}

\begin{theorem}[Theorems~1.75 and 1.77 of \cite{B:P:01}]
If $a,b,c \in \mathbb{T}$, $a
\le c \le b$, $\alpha \in \mathbb{R}$, and $f,g \in
C_{\textrm{rd}}(\mathbb{T}, \mathbb{R})$, then
\begin{enumerate}
\item $\displaystyle \int_{a}^{b}\left(f(t) + g(t) \right)
    \Delta t= \int_{a}^{b}f(t)\Delta t +
    \int_{a}^{b}g(t)\Delta t$;

\item $\displaystyle \int_{a}^{b} \alpha f(t)\Delta t =\alpha
    \int_{a}^{b}f(t)\Delta t$;

\item $\displaystyle \int_{a}^{b}  f(t)\Delta t = -
    \int_{b}^{a} f(t)\Delta t$;

\item $\displaystyle \int_{a}^{a}  f(t)\Delta t=0$;

\item $\displaystyle \int_{a}^{b}  f(t)\Delta t =
    \int_{a}^{c}  f(t)\Delta t + \int_{c}^{b} f(t)\Delta t$;

\item If $f(t)> 0$ for all $a \leq  t< b$, then $
    \displaystyle \int_{a}^{b}  f(t)\Delta t > 0$;

\item If $t \in \mathbb{T}^\kappa$, then $\displaystyle
\int_{t}^{\sigma(t)}  f(\tau)\Delta \tau =\mu(t)f(t)$.

\end{enumerate}
\end{theorem}

We have presented only the very basic concepts of the theory of time scales.
Nowadays the time scales theory covers nonlinear and higher order dynamic equations,
boundary value problems, dynamic inequalities, symplectic dynamical systems, etc.
An analogous theory has been later developed for the ``nabla derivative'', denoted $f^\nabla$,
which is a generalization of the backward difference operator from discrete calculus
\cite{A:G:02,B:P:01}. This is the subject of our next section.
The nabla calculus seems to be particularly useful
as a modeling technique in the calculus of variations
with applications to economics \cite{A:B:L:06,A:U:08,MyID:170}.


\subsection{The nabla approach}

In order to introduce the definition of nabla derivative, we define
a new set $\mathbb{T}_\kappa$ which is derived from $\mathbb{T}$ as
follows: if  $\mathbb{T}$ has a right-scattered minimum $m$, then
$\mathbb{T}_\kappa=\mathbb{T}\setminus\{m\}$; otherwise,
$\mathbb{T}_\kappa= \mathbb{T}$. In order to simplify expressions,
and similarly as done with composition with $\sigma$, we define
$f^{\rho}(t) := f(\rho(t))$.

\begin{definition}
We say that a function $f:\mathbb{T}\rightarrow\mathbb{R}$ is
\emph{nabla differentiable} at $t\in\mathbb{T}_\kappa$ if there is a
number $f^{\nabla}(t)$ such that for all $\varepsilon>0$ there
exists a neighborhood $U$ of $t$ such that
$$
|f^\rho(t)-f(s)-f^{\nabla}(t)(\rho(t)-s)|
\leq\varepsilon|\rho(t)-s|
$$
for all $s\in U$. We call $f^{\nabla}(t)$
the \emph{nabla derivative} of $f$ at $t$. Moreover,
we say that $f$ is \emph{nabla differentiable} on $\mathbb{T}$
provided $f^{\nabla}(t)$ exists for all $t \in \mathbb{T}_\kappa$.
\end{definition}

\begin{theorem}(Theorem~8.39 in \cite{B:P:01})
\label{propriedades derivada} Let $\mathbb{T}$ be a time scale,
$f:\mathbb{T}\rightarrow\mathbb{R}$, and $t\in\mathbb{T}_\kappa$.
The following holds:
\begin{enumerate}
\item If $f$ is nabla differentiable at $t$, then $f$ is
    continuous at $t$.
\item If $f$ is continuous at $t$ and $t$ is left-scattered,
    then $f$ is nabla differentiable at $t$ and
$$f^{\nabla}(t)=\frac{f(t)-f(\rho(t))}{t-\rho(t)}.$$

\item If $t$ is left-dense, then $f$ is nabla differentiable
    at $t$ if and only if the limit
$$
\lim_{s\rightarrow t} \frac{f(t)-f(s)}{t-s}
$$
exists as a finite number. In this case,
$$
f^\nabla(t)=\lim_{s\rightarrow t} \frac{f(t)-f(s)}{t-s}.
$$
\item If $f$ is nabla differentiable at $t$, then
$$
f^\rho(t)=f(t)-\nu(t)f^\nabla(t).
$$
\end{enumerate}
\end{theorem}

\begin{remark}
When $\mathbb{T}=\mathbb{R}$, then $f:\mathbb{R} \rightarrow
\mathbb{R}$ is nabla differentiable  at $t \in \mathbb{R}$ if and
only if $$ \displaystyle f^{\nabla}(t)= \lim_{s\rightarrow
t}\frac{f(t)-f(s)}{t-s}$$ exists, \textrm{i.e.}, if and only if $f$
is differentiable at $t$ in the ordinary sense. When
$\mathbb{T}=\mathbb{Z}$, then $f:\mathbb{Z} \rightarrow \mathbb{R}$
is always nabla differentiable at $t \in \mathbb{Z}$ and
$$
f^{\nabla}(t)=\frac{f(t)-f(\rho(t))}{t-\rho(t)}= f(t)-f(t-1)=:\nabla
f(t),
$$
\textrm{i.e.}, $\nabla$ is the usual backward difference
operator defined by the last equation above. For any time scale
$\mathbb{T}$, when $f$ is a constant, then $f^{\nabla}=0$; if
$f(t)=k t$ for some constant $k$, then $f^{\nabla}=k$.
\end{remark}

\begin{theorem}(Theorem~8.41 in \cite{B:P:01})
Suppose $f,g:\mathbb{T}\rightarrow\mathbb{R}$ are nabla
differentiable at $t\in\mathbb{T}_\kappa$. Then,
\begin{enumerate}
\item the sum $f+g:\mathbb{T}\rightarrow\mathbb{R}$ is nabla
differentiable at $t$ and
$$(f+g)^{\nabla}(t)=f^{\nabla}(t) + g^{\nabla}(t) \, ;$$

\item for any constant $\alpha$, $\alpha
f:\mathbb{T}\rightarrow\mathbb{R}$ is nabla differentiable at $t$
and $$(\alpha f)^{\nabla} (t)=\alpha f^{\nabla}(t)\, ;$$

\item the product $fg:\mathbb{T}\rightarrow\mathbb{R}$ is
    nabla differentiable at $t$ and
$$
\begin{array}{rcl}
(fg)^{\nabla}(t)& = & f^{\nabla}(t)g(t) + f^{\rho}(t)g^{\nabla}(t)\\
&&\\ &=& f^{\nabla}(t)g^{\rho}(t)+ f(t)g^{\nabla}(t).
\end{array}
$$
\end{enumerate}
\end{theorem}

\begin{definition}
A function $F:\mathbb{T}\rightarrow\mathbb{R}$ is called a
\emph{nabla antiderivative} of $f:\mathbb{T}\rightarrow\mathbb{R}$
provided $F^{\nabla}(t)=f(t)$ for all $t \in \mathbb{T}_\kappa$.
In this case we define the \emph{nabla integral} of $f$ from $a$ to
$b$ ($a,b \in \mathbb{T}$) by
$$
\int_{a}^{b}f(t)\nabla t:=F(b)-F(a) \, .
$$
\end{definition}

In order to exhibit a class of functions that possess a nabla
antiderivative, the following definition is introduced.

\begin{definition}
Let $\mathbb{T}$ be a time scale,
$f:\mathbb{T}\rightarrow\mathbb{R}$. We say that function $f$ is
\emph{ld-continuous} if it is continuous at left-dense points
and its right-sided limits exist (finite) at all right-dense points.
\end{definition}

\begin{theorem}(Theorem~8.45 in \cite{B:P:01})
Every ld-continuous function has a nabla
antiderivative. In particular, if $a \in \mathbb{T}$, then the
function $F$ defined by
$$
F(t)= \int_{a}^{t}f(\tau)\nabla\tau, \quad t \in \mathbb{T} \, ,
$$
is a nabla antiderivative of $f$.
\end{theorem}

The set of all ld-continuous functions
$f:\mathbb{T}\rightarrow\mathbb{R}$ is denoted by
$C^0_{\textrm{ld}} = C_{\textrm{ld}}
= C_{\textrm{ld}}(\mathbb{T})
= C_{\textrm{ld}}(\mathbb{T}, \mathbb{R})$,
and the set of all nabla differentiable functions
with ld-continuous derivative by
$C_{\textrm{ld}}^1 = C_{\textrm{ld}}^1(\mathbb{T}, \mathbb{R})$.
In general, we say that $f \in C_{\textrm{ld}}^n$, $n \in \mathbb{N}$,
if $f^\Delta \in C_{\textrm{ld}}^{n-1}$.

\begin{theorem}(Theorems~8.46 and 8.47 in \cite{B:P:01})
If $a$, $b$, $c \in \mathbb{T}$, $a \le c \le b$,
$\alpha \in \mathbb{R}$,
and $f,g \in C_{\textrm{ld}}(\mathbb{T}, \mathbb{R})$, then
\begin{enumerate}
\item $\displaystyle \int_{a}^{b}\left(f(t) + g(t) \right)
    \nabla t= \int_{a}^{b}f(t)\nabla t +
    \int_{a}^{b}g(t)\nabla t$;

\item $\displaystyle \int_{a}^{b} \alpha f(t)\nabla t =\alpha
    \int_{a}^{b}f(t)\nabla t$;

\item $\displaystyle \int_{a}^{b}  f(t)\nabla t = -
    \int_{b}^{a} f(t)\nabla t$;

\item $\displaystyle \int_{a}^{a}  f(t)\nabla t=0$;

\item $\displaystyle \int_{a}^{b}  f(t)\nabla t =
    \int_{a}^{c}  f(t)\nabla t + \int_{c}^{b} f(t)\nabla t$;

\item If $f(t)> 0$ for all $a < t \leq b$, then $
    \displaystyle \int_{a}^{b}  f(t)\nabla t > 0$;

\item If $t \in \mathbb{T}_\kappa$, then
$\displaystyle \int_{\rho(t)}^{t}f(\tau)\nabla\tau=\nu(t)f(t)$.

\end{enumerate}
\end{theorem}

\begin{example}
Let $a, b \in \mathbb{T}$ and $f\in C_{\textrm{ld}}(\mathbb{T}, \mathbb{R})$.
For $\mathbb{T}=\mathbb{R}$, then $\int_{a}^{b}f(t)\nabla t =
\int_{a}^{b}f(t) dt$, where the integral on the right side is the
usual Riemann integral. For $\mathbb{T}=\mathbb{Z}$, then
$\displaystyle \int_{a}^{b}f(t)\nabla t = \sum_{t=a+1}^{b}f(t)$ if
$a<b$, $\displaystyle \int_{a}^{b}f(t)\nabla t=0$ if $a=b$, and
$\displaystyle \int_{a}^{b}f(t)\nabla t = - \sum_{t=b+1}^{a}f(t)$ if
$a>b$.
\end{example}

For more on the nabla calculus we refer the reader to
\cite[Chap.~3]{B:P:03}.


\subsection{Preliminaries to variational calculus}

Our goal is to obtain Euler--Lagrange type equations. Similar
to the classical calculus of variations \cite{GelfandFomin},
integration by parts will play an important role.
If functions $f,g : \mathbb{T}\rightarrow\mathbb{R}$
are delta and nabla differentiable with continuous derivatives,
then the following formulas of integration by parts hold \cite{B:P:01}:
\begin{equation}
\label{intBP}
\begin{split}
\int_{a}^{b}f^\sigma(t) g^{\Delta}(t)\Delta t
&=\left.(fg)(t)\right|_{t=a}^{t=b}
-\int_{a}^{b}f^{\Delta}(t)g(t)\Delta t \, , \\
\int_{a}^{b}f(t)g^{\Delta}(t)\Delta t
&=\left.(fg)(t)\right|_{t=a}^{t=b}
-\int_{a}^{b}f^{\Delta}(t)g^\sigma(t)\Delta t \, , \\
\int_{a}^{b}f^\rho(t)g^{\nabla}(t)\nabla t
&=\left.(fg)(t)\right|_{t=a}^{t=b}
-\int_{a}^{b}f^{\nabla}(t)g(t)\nabla t \, ,\\
\int_{a}^{b}f(t)g^{\nabla}(t)\nabla t
&=\left.(fg)(t)\right|_{t=a}^{t=b}
-\int_{a}^{b}f^{\nabla}(t)g^\rho(t)\nabla t \, .
\end{split}
\end{equation}

The following fundamental lemma of the calculus of variations
on time scales, involving a nabla derivative
and a nabla integral, has been proved in \cite{NM:T}.

\begin{lemma}{\rm (The nabla Dubois--Reymond lemma \cite[Lemma~14]{NM:T}).}
\label{DBRL:n}
Let $f \in C_{\textrm{ld}}([a,b], \mathbb{R})$. If
$$
\int_{a}^{b} f(t)\eta^{\nabla}(t)\nabla t=0
$$
for all $\eta \in C_{\textrm{ld}}^1([a,b],
\mathbb{R})$ with $\eta(a)=\eta(b)=0$,
then $f(t) = c$, $t\in [a,b]_\kappa$,
for some constant $c$.
\end{lemma}

Lemma~\ref{DBRL:d} is the analogous
delta version of Lemma~\ref{DBRL:n}.

\begin{lemma}{\rm (The delta Dubois--Reymond lemma \cite{B:04}).}
\label{DBRL:d}
Let $g\in C_{\textrm{rd}}([a,b], \mathbb{R})$. If
$$
\int_{a}^{b} g(t) \eta^\Delta(t)\Delta t=0
$$
for all $\eta \in C_{\textrm{rd}}^1$ with
$\eta(a)=\eta(b)=0$, then $g(t)=c$
on $[a,b]^\kappa$ for some $c\in\mathbb{R}$.
\end{lemma}

Proposition~\ref{prop:rel:der}
gives a relationship between delta
and nabla derivatives.

\begin{proposition}{\rm (Theorems~2.5 and 2.6 of \cite{A:G:02}).}
\label{prop:rel:der}
(i) If $f : \mathbb{T} \rightarrow \mathbb{R}$ is delta differentiable
on $\mathbb{T}^\kappa$ and $f^\Delta$ is continuous on $\mathbb{T}^\kappa$,
then $f$ is nabla differentiable on $\mathbb{T}_\kappa$ and
\begin{equation}
\label{eq:chgN_to_D}
f^\nabla(t)=\left(f^\Delta\right)^\rho(t) \quad \text{for all }
t \in \mathbb{T}_\kappa \, .
\end{equation}
(ii) If $f : \mathbb{T} \rightarrow \mathbb{R}$ is nabla differentiable
on $\mathbb{T}_\kappa$ and $f^\nabla$ is continuous on $\mathbb{T}_\kappa$,
then $f$ is delta differentiable on $\mathbb{T}^\kappa$ and
\begin{equation}
\label{eq:chgD_to_N}
f^\Delta(t)=\left(f^\nabla\right)^\sigma(t) \quad \text{for all }
t \in \mathbb{T}^\kappa \, .
\end{equation}
\end{proposition}

\begin{remark}
Note that, in general, $f^\nabla(t) \ne f^\Delta\left(\rho(t)\right)$
and $f^\Delta(t) \ne f^\nabla\left(\sigma(t)\right)$.
In Proposition~\ref{prop:rel:der} the assumptions on the
continuity of $f^\Delta$ and $f^\nabla$ are crucial.
\end{remark}

\begin{proposition}{\rm (\cite[Theorem~2.8]{A:G:02}).}
\label{eq:prop}
Let $a, b \in \mathbb{T}$ with $a \le b$ and let $f$
be a continuous function on $[a, b]$. Then,
\begin{equation*}
\begin{split}
\int_a^b f(t)\Delta t
&= \int_a^{\rho(b)} f(t)\Delta t
+ (b - \rho(b))f^\rho(b) \, , \\
\int_a^b f(t)\Delta t
&= (\sigma(a) - a) f(a)
+ \int_{\sigma(a)}^b f(t)\Delta t \, , \\
\int_a^b f(t)\nabla t
&= \int_a^{\rho(b)} f(t)\nabla t
+ (b - \rho(b)) f(b) \, , \\
\int_a^b f(t)\nabla t
&= (\sigma(a) - a) f^\sigma(a)
+ \int_{\sigma(a)}^b f(t)\nabla t \, .
\end{split}
\end{equation*}
\end{proposition}

We end our brief review of the calculus on time scales
with a relationship between the delta and nabla integrals.

\begin{proposition}{\rm (\cite[Proposition~7]{G:G:S:05}).}
If function $f : \mathbb{T} \rightarrow \mathbb{R}$
is continuous, then for all $a, b \in \mathbb{T}$
with $a < b$ we have
\begin{gather}
\int_a^b f(t) \Delta t = \int_a^b f^\rho(t) \nabla t \, , \label{eq:DtoN}\\
\int_a^b f(t) \nabla t = \int_a^b f^\sigma(t) \Delta t \, . \label{eq:NtoD}
\end{gather}
\end{proposition}


\section{The $\Delta$-$\nabla$ calculus of variations}

We consider the problem of extremizing
the variational functional \eqref{eq:P}
subject to given boundary conditions
$y(a) = \alpha$ and $y(b) = \beta$:
\begin{equation}
\label{problem:P}
\begin{gathered}
\mathcal{J}(y) =
\left(\int_a^b L_{\Delta}[y](t) \Delta t\right)
\left(\int_a^b L_{\nabla}\{y\}(t) \nabla t\right) \rightarrow
\textrm{extr} \\
y \in C_{\diamond}^1\left([a,b],\mathbb{R}\right) \\
y(a) = \alpha \, , \quad y(b) = \beta \, ,
\end{gathered}
\end{equation}
where $C_{\diamond}^1\left([a,b],\mathbb{R}\right)$
denote the class of functions
$y : [a,b]\rightarrow\mathbb{R}$  with
$y^\Delta$ continuous on $[a,b]^\kappa$
and $y^\nabla$ continuous on $[a,b]_\kappa$.

Before presenting the Euler--Lagrange equations for
problem \eqref{problem:P} we introduce the definition
of weak local extremum.

\begin{definition} We say that $\hat{y}\in C_{\diamond}^{1}([a,b],
\mathbb{R})$ is a weak local minimizer (respectively weak local
maximizer) for problem \eqref{problem:P} if there exists
$\delta >0$ such that
$$
\mathcal{J}(\hat{y})\leq \mathcal{J}(y) \quad
(\text{respectively} \ \   \mathcal{J}(\hat{y})
\geq \mathcal{J}(y))
$$
for all $y \in C_{\diamond}^{1}([a,b], \mathbb{R})$
satisfying the boundary
conditions $y(a) = \alpha$ and $y(b) = \beta$, and
$$
\parallel y - \hat{y}\parallel_{1,\infty} < \delta \, ,
$$
where
$$
\parallel y\parallel_{1,\infty}:=
\parallel y^{\sigma}\parallel_{\infty}
+ \parallel y^{\rho}\parallel_{\infty}
+ \parallel y^{\Delta}\parallel_{\infty}
+ \parallel y^{\nabla}\parallel_{\infty}
$$
and
$$ \parallel y\parallel_{\infty} :=\sup_{t \in
[a,b]_{\kappa}^{\kappa}}\mid y(t) \mid.
$$
\end{definition}


\subsection{Euler--Lagrange equations}
\label{subsec:mr:el}

Theorem~\ref{thm:mr} gives two different forms
for the Euler--Lagrange equation on time scales
associated with the variational problem \eqref{problem:P}.

\begin{theorem}{\rm (The general Euler--Lagrange equations
on time scales).}
\label{thm:mr}
If $\hat{y} \in C_{\diamond}^1$ is a weak local extremizer of problem
\eqref{problem:P}, then $\hat{y}$ satisfies
the following delta-nabla integral equations:
\begin{multline}
\label{eq:EL1}
\mathcal{J}_\nabla(\hat{y})
\left(\partial_3 L_\Delta[\hat{y}](\rho(t))
-\int_{a}^{\rho(t)} \partial_2 L_\Delta[\hat{y}](\tau) \Delta\tau\right)\\
+
\mathcal{J}_\Delta(\hat{y})
\left(\partial_3 L_\nabla\{\hat{y}\}(t)
-\int_{a}^{t} \partial_2 L_\nabla\{\hat{y}\}(\tau) \nabla\tau\right)
= \text{const}
\end{multline}
for all $t \in [a,b]_\kappa$;
\begin{multline}
\label{eq:EL2}
\mathcal{J}_\Delta(\hat{y})
\left(\partial_3 L_\nabla\{\hat{y}\}(\sigma(t))
-\int_{a}^{\sigma(t)} \partial_2 L_\nabla\{\hat{y}\}(\tau) \nabla\tau\right)\\
+
\mathcal{J}_\nabla(\hat{y})
\left(\partial_3 L_\Delta[\hat{y}](t)
-\int_{a}^{t} \partial_2 L_\Delta[\hat{y}](\tau) \Delta\tau\right)
= \text{const}
\end{multline}
for all $t \in [a,b]^\kappa$.
\end{theorem}

\begin{remark}
In the classical context (\textrm{i.e.},
when $\mathbb{T} = \mathbb{R}$) the
necessary conditions \eqref{eq:EL1}
and \eqref{eq:EL2} coincide with the Euler--Lagrange
equations \eqref{eq:EL:P} given in \cite{Pedregal}.
\end{remark}

\begin{proof}
Suppose that $\mathcal{J}$
has a weak local extremum at $\hat{y}$. We
consider the value of $\mathcal{J}$ at nearby
functions $\hat{y} + \varepsilon \eta$,
where $\varepsilon\in \mathbb{R}$ is a small parameter,
$\eta \in C_{\diamond}^{1}([a,b],\mathbb{R})$ with $\eta(a)=\eta(b)=0$.
Thus, function $\phi(\varepsilon)
= \mathcal{J}(\hat{y} + \varepsilon \eta)$
has an extremum at $\varepsilon = 0$. Using the first-order necessary
optimality condition $\left.\phi'(\varepsilon)\right|_{\varepsilon = 0} = 0$
we obtain:
\begin{multline}
\label{eq:prf:+}
\hspace*{-0.5cm}
0=\mathcal{J}_\Delta(\hat{y}) \int_a^b
\left(\partial_2 L_\nabla\{\hat{y}\}(t) \eta^\rho(t)
+ \partial_3 L_\nabla\{\hat{y}\}(t) \eta^\nabla(t)\right)\nabla t\\
+ \mathcal{J}_\nabla(\hat{y}) \int_a^b
\left(\partial_2 L_\Delta[\hat{y}](t) \eta^\sigma(t)
+ \partial_3 L_\Delta[\hat{y}](t) \eta^\Delta(t)\right) \Delta t.
\end{multline}
Let
\begin{equation*}
A(t) = \int_a^t \partial_2 L_\Delta[\hat{y}](\tau) \Delta\tau,
\quad
B(t) = \int_a^t \partial_2 L_\nabla\{\hat{y}\}(\tau) \nabla\tau.
\end{equation*}
Then, $A^\Delta(t) = \partial_2 L_\Delta[\hat{y}](t)$,
$B^\nabla(t) = \partial_2 L_\nabla\{\hat{y}\}(t)$,
and the first and third integration by parts formula
in \eqref{intBP} tell us, respectively, that
\begin{equation*}
\begin{split}
\int_a^b &\partial_2 L_\Delta[\hat{y}](t) \eta^\sigma(t) \Delta t
= \int_a^b A^\Delta(t) \eta^\sigma(t) \Delta t\\
&= \left. A(t) \eta(t)\right|_{t=a}^{t=b} - \int_a^b A(t) \eta^\Delta(t) \Delta t\\
&= - \int_a^b A(t) \eta^\Delta(t) \Delta t
\end{split}
\end{equation*}
and
\begin{equation*}
\begin{split}
\int_a^b &\partial_2 L_\nabla\{\hat{y}\}(t) \eta^\rho(t) \nabla t
= \int_a^b B^\nabla(t) \eta^\rho(t) \nabla t\\
&= \left. B(t) \eta(t)\right|_{t=a}^{t=b}
- \int_a^b B(t) \eta^\nabla(t) \nabla t\\
&= - \int_a^b B(t) \eta^\nabla(t) \nabla t \, .
\end{split}
\end{equation*}
If we denote $f(t) = \partial_3 L_\Delta[\hat{y}](t) - A(t)$
and $g(t) = \partial_3 L_\nabla\{\hat{y}\}(t) - B(t)$,
then we can write the necessary optimality condition
\eqref{eq:prf:+} in the form
\begin{equation}
\label{eq:prf:+:aftIP}
\mathcal{J}_\nabla(\hat{y}) \int_a^b f(t) \eta^\Delta(t) \Delta t
+ \mathcal{J}_\Delta(\hat{y}) \int_a^b g(t) \eta^\nabla(t) \nabla t = 0 \, .
\end{equation}
We now split the proof in two parts:
(i) we prove \eqref{eq:EL1} transforming the delta integral
in \eqref{eq:prf:+:aftIP} to a nabla integral by means of
\eqref{eq:DtoN}; (ii) we prove \eqref{eq:EL2} transforming
the nabla integral in \eqref{eq:prf:+:aftIP}
to a delta integral by means of \eqref{eq:NtoD}.
(i) By \eqref{eq:DtoN} the necessary optimality
condition \eqref{eq:prf:+:aftIP} is equivalent to
\begin{equation*}
\int_a^b \left(\mathcal{J}_\nabla(\hat{y})
f^\rho(t) (\eta^\Delta)^\rho(t) + \mathcal{J}_\Delta(\hat{y})
g(t) \eta^\nabla(t)\right) \nabla t = 0
\end{equation*}
and by \eqref{eq:chgN_to_D} to
\begin{equation}
\label{eq:bef:FL1}
\int_a^b \left(\mathcal{J}_\nabla(\hat{y}) f^\rho(t)
+ \mathcal{J}_\Delta(\hat{y}) g(t)\right)
\eta^\nabla(t) \nabla t = 0 \, .
\end{equation}
Applying Lemma~\ref{DBRL:n} to \eqref{eq:bef:FL1}
we prove \eqref{eq:EL1}:
\begin{equation}
\label{eq:S12}
\mathcal{J}_\nabla(\hat{y}) f^\rho(t)
+ \mathcal{J}_\Delta(\hat{y}) g(t) = c \quad \forall t \in [a,b]_\kappa \, ,
\end{equation}
where $c$ is a constant.
(ii) By \eqref{eq:NtoD} the necessary optimality
condition \eqref{eq:prf:+:aftIP} is equivalent to
\begin{equation*}
\int_a^b \left(\mathcal{J}_\nabla(\hat{y}) f(t) \eta^\Delta(t)
+ \mathcal{J}_\Delta(\hat{y}) g^\sigma(t)
\left(\eta^\nabla\right)^\sigma(t)\right) \Delta t = 0
\end{equation*}
and by \eqref{eq:chgD_to_N} to
\begin{equation}
\label{eq:bef:FL2}
\int_a^b \left(\mathcal{J}_\nabla(\hat{y}) f(t)
+ \mathcal{J}_\Delta(\hat{y}) g^\sigma(t)\right)
\eta^\Delta(t) \Delta t = 0 \, .
\end{equation}
Applying Lemma~\ref{DBRL:d} to \eqref{eq:bef:FL2}
we prove \eqref{eq:EL2}:
\begin{equation}
\label{eq:S13}
\mathcal{J}_\nabla(\hat{y}) f(t)
+ \mathcal{J}_\Delta(\hat{y}) g^\sigma(t)
= c \quad \forall t \in [a,b]^\kappa \, ,
\end{equation}
where $c$ is a constant.
\end{proof}

\begin{corollary}
Let $L_\Delta\left(t,y^\sigma,y^\Delta\right) = L_\Delta(t)$
and $\mathcal{J}_\Delta(\hat{y}) \ne 0$
(this is true, \textrm{e.g.}, for
$L_\Delta \equiv \frac{1}{b-a}$ for which $\mathcal{J}_\Delta
= 1$; \textrm{cf.} Remark~\ref{obs}). Then, $\partial_2 L_\Delta
= \partial_3 L_\Delta = 0$ and the Euler--Lagrange equation \eqref{eq:EL1} takes the form
\begin{equation}
\label{cor:EL1}
\partial_3 L_\nabla\{\hat{y}\}(t)
-\int_{a}^{t} \partial_2 L_\nabla\{\hat{y}\}(\tau) \nabla\tau
= \text{const}
\end{equation}
for all $t \in [a,b]_\kappa$.
\end{corollary}

\begin{remark}
If $\hat{y} \in C_{\textrm{ld}}^2$, then
nabla-differentiating \eqref{cor:EL1} we obtain
the Euler--Lagrange differential equation
\eqref{eq:el:n} as proved in \cite{NM:T}:
\begin{equation*}
\frac{\nabla}{\nabla t}
\partial_3 L_\nabla\{\hat{y}\}(t)
- \partial_2 L_\nabla\{\hat{y}\}(t) = 0 \quad \forall t \in [a,b]_{\kappa^2} \, .
\end{equation*}
\end{remark}

\begin{corollary}
Let $L_\nabla\left(t,y^\rho,y^\nabla\right) = L_\nabla(t)$
and $\mathcal{J}_\nabla(\hat{y}) \ne 0$
(this is true, \textrm{e.g.}, for
$L_\nabla \equiv \frac{1}{b-a}$ for which $\mathcal{J}_\nabla
= 1$; \textrm{cf.} Remark~\ref{obs}). Then, $\partial_2 L_\nabla
= \partial_3 L_\nabla = 0$ and the
Euler--Lagrange equation \eqref{eq:EL2} takes the form
\begin{equation}
\label{cor:EL2}
\partial_3 L_\Delta[\hat{y}](t)
-\int_{a}^{t} \partial_2 L_\Delta[\hat{y}](\tau) \Delta\tau
= \text{const}
\end{equation}
for all $t \in [a,b]^\kappa$.
\end{corollary}

\begin{remark}
If $\hat{y} \in C_{\textrm{rd}}^2$,
then delta-differentiating \eqref{cor:EL2} we obtain the Euler--Lagrange
differential equation \eqref{eq:el:d} as proved in \cite{B:04}:
\begin{equation*}
\frac{\Delta}{\Delta t}
\partial_3 L_\Delta[\hat{y}](t)
- \partial_2 L_\Delta[\hat{y}](t) = 0 \quad \forall t \in [a,b]^{\kappa^2} \, .
\end{equation*}
\end{remark}

\begin{example}
\label{ex:first:simp:ex}
Let $\mathbb{T}$ be a time scale with $0$, $\xi \in \mathbb{T}$,
$0 < \xi$, and $\left(\mathbb{T}\setminus \{0,\xi\}\right)\cap [0,\xi] \ne \emptyset$.
Consider the problem
\begin{equation}
\label{ex:1}
\begin{gathered}
\mathcal{J}(y)=\left(\int_{0}^{\xi}(y^\Delta(t))^2\Delta t\right)
\left(\int_{0}^{\xi}\left(y^\nabla(t))^2\right)\nabla t\right) \rightarrow \min\\
y(0)=0, \quad y(\xi)=\xi \, .
\end{gathered}
\end{equation}
Since
\begin{equation*}
L_{\Delta}=(y^\Delta)^2,\quad L_{\nabla}=(y^\nabla)^2
\end{equation*}
we have
\begin{equation*}
\partial_2L_{\Delta}=0,\
\partial_3L_{\Delta}=2y^\Delta, \
\partial_2L_{\nabla}=0, \
\partial_3L_{\nabla}=2y^\nabla.
\end{equation*}
Using equation \eqref{eq:EL2} of Theorem~\ref{thm:mr}
we get the following delta-nabla differential equation:
\begin{equation}\label{ex:2}
2Ay^{\Delta}(t)+2By^{\nabla}(\sigma(t))=C,
\end{equation}
where $C\in\mathbb{R}$ and $A$, $B$ are the values of functionals
$\mathcal{J}_{\nabla}$ and $\mathcal{J}_{\Delta}$ in a solution of
problem \eqref{ex:1}, respectively. From \eqref{eq:chgD_to_N}
we can write equation \eqref{ex:2} in the form
\begin{equation}
\label{ex:3}
2Ay^{\Delta}(t)+2By^\Delta=C.
\end{equation}
Observe that $A+B$ cannot be equal to $0$. Thus, solving equation
\eqref{ex:3} subject to the boundary conditions $y(0)=0$ and
$y(\xi)=\xi$ we get $y(t)=t$ as a candidate local minimizer for the
problem \eqref{ex:1}.
\end{example}

\begin{example}
Consider the problem
\begin{equation}
\label{ex:product}
\begin{gathered}
\mathcal{J}(y)=\left(\int_{0}^{1}ty^{\Delta}(t) \Delta
     t\right)\left(\int_{0}^{1}(y^{\nabla}(t))^2\nabla
     t\right) \rightarrow \textrm{extr} \\
     y(0)=0, \quad y(1)=1.
\end{gathered}
\end{equation}
Since
\begin{equation*}
L_{\Delta}=ty^\Delta,\quad L_{\nabla}=(y^\nabla)^2
\end{equation*}
we have
\begin{equation*}
\partial_2L_{\Delta}=0,\quad \partial_3L_{\Delta}=t,
\quad \partial_2L_{\nabla}=0,
\quad \partial_3L_{\nabla}=2y^\nabla.
\end{equation*}
Using equation \eqref{eq:EL2} of Theorem~\ref{thm:mr}
and relation \eqref{eq:chgD_to_N}, we get
the following delta differential equation:
\begin{equation}
\label{ex:product:1}
At+2By^{\Delta}(t)=C,
\end{equation}
where $C\in\mathbb{R}$ and $A$, $B$ are values of the functionals
$\mathcal{J}_{\nabla}$ and $\mathcal{J}_{\Delta}$ in a solution of
\eqref{ex:product}, respectively. Observe that $A\neq0$, so that $B$ is
also nonzero. A solution of \eqref{ex:product:1} depends on the
time scale. Let us solve, for example, this equation on
$\mathbb{T}=\mathbb{R}$ and on
$\mathbb{T}=\left\{0,\frac{1}{2},1\right\}$. On
$\mathbb{T}=\mathbb{R}$ we obtain
\begin{equation}\label{sol:P}
y(t)=-\frac{A}{4B}t^2+\frac{4B+A}{4B}t.
\end{equation}
Substituting \eqref{sol:P} into functionals $\mathcal{J}_{\nabla}$
and $\mathcal{J}_{\Delta}$ gives
\begin{equation}\label{equation:A,B}
\begin{cases}
\frac{48B^2+A^2}{48B^2}=A\\
\frac{12B-A}{24B}=B.
\end{cases}
\end{equation}
Solving the system of equations \eqref{equation:A,B} we obtain
\begin{gather*}
\begin{cases}
A=0\\
B=0,
\end{cases}
\quad
\begin{cases}
A=\frac{4}{3}\\
B=\frac{1}{3}.
\end{cases}
\end{gather*}
Therefore,
\begin{equation*}
y(t)=-t^2+2t
\end{equation*}
is a candidate extremizer for problem \eqref{ex:product} on
$\mathbb{T}=\mathbb{R}$. Note that nothing can be concluded
from Theorem~\ref{thm:mr} as to whether $y$ gives a minimum,
a maximum, or neither of these, for $\mathcal{J}$.

The solution of \eqref{ex:product:1} on
$\mathbb{T}=\left\{0,\frac{1}{2},1\right\}$ is
\begin{gather}\label{ex:product:2}
y(t)=
\begin{cases}
0 & \text{ if } t=0 \\
 \frac{1}{2}+\frac{A}{16B} & \text{ if } t=\frac{1}{2}\\
 1 & \text{ if } t=1.
\end{cases}
\end{gather}
Constants $A$ and $B$ are determined by substituting
\eqref{ex:product:2} into functionals $\mathcal{J}_{\nabla}$ and
$\mathcal{J}_{\Delta}$. The resulting system of equations is
\begin{equation}\label{equation:T}
    \begin{cases}
1+\frac{A^2}{64B^2}=A\\
\frac{1}{4}-\frac{A}{32B}=B.
\end{cases}
\end{equation}
Since system of equations \eqref{equation:T} has no real solutions,
we conclude that there exists no extremizer for problem
\eqref{ex:product} on $\mathbb{T}=\left\{0,\frac{1}{2},1\right\}$
among the set of functions that we consider to be admissible.
\end{example}


\subsection{Natural boundary conditions}
\label{subsec:mr:nbc}

We consider now the situations when we want to minimize
or maximize the variational functional $\mathcal{J}$ but
$y(a)$ and/or $y(b)$ are free. Then,
$$
\left. A(t) \eta(t)\right|_{t=a}^{t=b}
= A(b) \eta(b) - A(a) \eta(a) = A(b) \eta(b)
$$
and
$$
\left. B(t) \eta(t)\right|_{t=a}^{t=b}
= B(b) \eta(b) - B(a) \eta(a) = B(b) \eta(b)
$$
do not vanish necessarily, and conditions
\eqref{eq:bef:FL1} and \eqref{eq:bef:FL2} take,
respectively, the form
\begin{multline}
\label{eq:new16}
\left(A(b) + B(b)\right) \eta(b)\\
+ \int_a^b \left(\mathcal{J}_\nabla(\hat{y}) f^\rho(t)
+ \mathcal{J}_\Delta(\hat{y}) g(t)\right)
\eta^\nabla(t) \nabla t = 0
\end{multline}
and
\begin{multline}
\label{eq:new17}
\left(A(b) + B(b)\right) \eta(b)\\
+ \int_a^b \left(\mathcal{J}_\nabla(\hat{y}) f(t)
+ \mathcal{J}_\Delta(\hat{y}) g^\sigma(t)\right)
\eta^\Delta(t) \Delta t = 0 \, .
\end{multline}
Since \eqref{eq:new16} and \eqref{eq:new17} are valid
for an arbitrary function $\eta$, in particular they hold
for the subclass of functions $\eta$ for which
$\eta(a) = \eta(b) = 0$. Thus, the same Euler--Lagrange
conditions \eqref{eq:S12} and \eqref{eq:S13}
of Theorem~\ref{thm:mr} are obtained and both
\eqref{eq:new16} and \eqref{eq:new17} simplify to
$$
\left(c + A(b) + B(b)\right) \eta(b) - c \eta(a) = 0 \, .
$$
If $y(a)$ is free, $\eta(a)$ is arbitrary and we conclude that
$c = 0$; if $y(b)$ is free, $\eta(b)$ is arbitrary
and we conclude that $c + A(b) + B(b)= 0$. The value of the constant
$c$ is obtained from \eqref{eq:S12} and \eqref{eq:S13}.
We remark that
\eqref{eq:S12} does not hold in $a$ in the case $a$
is right-scattered; and \eqref{eq:S13} does not hold in $b$
if $b$ is left-scattered. By \eqref{eq:S13} one has
$c = \mathcal{J}_\nabla(\hat{y}) f(a)
+ \mathcal{J}_\Delta(\hat{y}) g(\sigma(a))$,
and by \eqref{eq:S12}
$c = \mathcal{J}_\nabla(\hat{y}) f(\rho(b))
+ \mathcal{J}_\Delta(\hat{y}) g(b)$.
We have just proved the following result:

\begin{theorem}{\rm (The general natural boundary conditions
on time scales).}
\label{thm:mr2}
If $\hat{y}$ is a weak local extremizer of the variational
functional \eqref{eq:P}, then the Euler--Lagrange equations
\eqref{eq:EL1} and \eqref{eq:EL2} hold,
together with the natural condition
\begin{multline*}
\mathcal{J}_\Delta(\hat{y})
\left(\partial_3 L_\nabla\{\hat{y}\}(\sigma(a))
-\int_{a}^{\sigma(a)} \partial_2 L_\nabla\{\hat{y}\}(\tau) \nabla\tau\right)\\
+ \mathcal{J}_\nabla(\hat{y}) \partial_3 L_\Delta[\hat{y}](a) = 0
\end{multline*}
when $y(a)$ is free; together with the natural condition
\begin{multline}
\label{NBC:b}
\mathcal{J}_\nabla(\hat{y})
\left(\partial_3 L_\Delta[\hat{y}](\rho(b))
-\int_{a}^{\rho(b)} \partial_2 L_\Delta[\hat{y}](\tau) \Delta\tau\right)\\
+ \int_{a}^{b} \partial_2 L_\Delta[\hat{y}](\tau) \Delta\tau \\
+ \mathcal{J}_\Delta(\hat{y})
\left(\partial_3 L_\nabla\{\hat{y}\}(b)
-\int_{a}^{b} \partial_2 L_\nabla\{\hat{y}\}(\tau) \nabla\tau\right)\\
+ \int_{a}^{b} \partial_2 L_\nabla\{\hat{y}\}(\tau) \nabla\tau = 0
\end{multline}
when $y(b)$ is free.
\end{theorem}

Using Proposition~\ref{eq:prop},
the next corollary is obtained:

\begin{corollary}{\rm (The delta natural boundary conditions).}
If $\hat{y}$ is a weak local extremizer of the delta variational
functional \eqref{eq:Pd}, then the Euler--Lagrange equation
\eqref{cor:EL2} holds together with the natural condition
\begin{equation}
\label{eq:DNBCa}
\partial_3 L_\Delta[\hat{y}](a) = 0
\end{equation}
when $y(a)$ is free; together with the natural condition
\begin{equation}
\label{eq:DNBCb}
\partial_3 L_\Delta[\hat{y}](\rho(b))
+ \left(b - \rho(b)\right)
\partial_2 L_\Delta[\hat{y}](\rho(b)) = 0
\end{equation}
when $y(b)$ is free.
\end{corollary}

Analogous nabla natural conditions are also trivially
obtained from Theorem~\ref{thm:mr2}:

\begin{corollary}{\rm (The nabla natural boundary conditions).}
If $\hat{y}$ is a weak local extremizer of the nabla variational
functional \eqref{eq:Pn}, then the Euler--Lagrange equation
\eqref{cor:EL1} holds together with the natural condition
\begin{equation}
\label{eq:NNBCa}
\partial_3 L_\nabla\{\hat{y}\}(\sigma(a))
- \int_{a}^{\sigma(a)}
\partial_2 L_\nabla\{\hat{y}\}(\tau) \nabla\tau = 0
\end{equation}
when $y(a)$ is free; together with the natural condition
\begin{equation}
\label{eq:NNBCb}
\partial_3 L_\nabla\{\hat{y}\}(b) = 0
\end{equation}
when $y(b)$ is free.
\end{corollary}

\begin{remark}
The natural boundary condition \eqref{eq:NNBCa} can be written
in the equivalent form
$$
\partial_3 L_\nabla\{\hat{y}\}(\sigma(a))
- \left(\sigma(a) - a\right)
\partial_2 L_\nabla\{\hat{y}\}(a) = 0 \, .
$$
\end{remark}

In the classical context of the calculus of variations,
\textrm{i.e.}, when $\mathbb{T} = \mathbb{R}$, both
\eqref{eq:DNBCa} and \eqref{eq:NNBCa} reduce to
$\partial_3 L(a,y(a),y'(a)) = 0$ and both
\eqref{eq:DNBCb} and \eqref{eq:NNBCb} reduce to
$\partial_3 L(b,y(b),y'(b)) = 0$,
for some given Lagrangian $L$, which are the standard
natural boundary conditions of the calculus of variations
(\textrm{cf.}, \textrm{e.g.}, \cite{GelfandFomin}).
For the particular case $\mathbb{T} = \mathbb{R}$
we obtain from our Theorem~\ref{thm:mr2} a result
in \cite{Pedregal} that generalizes the classical
natural boundary conditions to functionals given
by the product of two integrals:

\begin{corollary}{\rm (\textrm{cf.} \cite{Pedregal}).}
If $\hat{y}$ is a weak extremizer
of the variational functional
\begin{equation*}
\begin{split}
\mathcal{I}(y) &= \mathcal{I}_1(y) \mathcal{I}_2(y)\\
&= \left(\int_a^b L_1\left(t,y(t),y'(t)\right) dt\right)\\
&\qquad \cdot \left(\int_a^b L_2\left(t,y(t),y'(t)\right) dt \right) \, ,
\end{split}
\end{equation*}
then the Euler--Lagrange equation
\begin{multline}
\label{eq:EL:P}
\mathcal{I}_2(\hat{y})
\left(\partial_3 L_1\left(t,\hat{y}(t),\hat{y}'(t)\right)
- \int_a^t \partial_2 L_1\left(t,\hat{y}(t),\hat{y}'(t)\right)
dt\right) \\
+ \mathcal{I}_1(\hat{y})
\left(\partial_3 L_2\left(t,\hat{y}(t),\hat{y}'(t)\right)
- \int_a^t \partial_2 L_2\left(t,\hat{y}(t),\hat{y}'(t)\right) dt \right)\\
= \text{const}
\end{multline}
holds for all $t \in [a,b]$. Moreover, the natural condition
\begin{multline*}
\mathcal{I}_2(\hat{y})
\partial_3 L_1\left(a,\hat{y}(a),\hat{y}'(a)\right)\\
+ \mathcal{I}_1(\hat{y})
\partial_3 L_2\left(a,\hat{y}(a),\hat{y}'(a)\right) = 0
\end{multline*}
holds when $y(a)$ is free; the natural condition
\begin{multline*}
\mathcal{I}_2(\hat{y})
\partial_3 L_1\left(b,\hat{y}(b),\hat{y}'(b)\right)
+ \mathcal{I}_1(\hat{y})
\partial_3 L_2\left(b,\hat{y}(b),\hat{y}'(b)\right)\\
+ \left(1 - \mathcal{I}_2(\hat{y})\right)
\int_a^b \partial_2 L_1\left(\tau,\hat{y}(\tau),\hat{y}'(\tau)\right) d\tau\\
+ \left(1 - \mathcal{I}_1(\hat{y})\right)
\int_a^b \partial_2 L_2\left(\tau,\hat{y}(\tau),\hat{y}'(\tau)\right) d\tau
= 0
\end{multline*}
holds when $y(b)$ is free.
\end{corollary}

\begin{example}
Let us consider the functional of Example~\ref{ex:first:simp:ex}
but where we are free to choose the value of $y$ at point $\xi$:
\begin{equation*}
\begin{gathered}
\mathcal{J}(y)=\left(\int_{0}^{\xi}(y^\Delta(t))^2\Delta t\right)
\left(\int_{0}^{\xi}\left(y^\nabla(t))^2\right)\nabla t\right) \rightarrow \min\\
y(0)=0, \quad y(\xi) \text{ free } \, .
\end{gathered}
\end{equation*}
Solving the Euler--Lagrange equation \eqref{ex:3}
with $y(0)=0$ gives the extremal
\begin{equation}
\label{eq:EXT}
\hat{y}(t) = \frac{C}{2(A+B)} t \, .
\end{equation}
In this case the natural boundary condition \eqref{NBC:b}
along \eqref{eq:EXT} simplifies to
$$
A\left(\frac{C}{A + B}\right) + B\left(\frac{C}{A+B}\right) = 0 \, ,
$$
that is, $C = 0$ and the minimum is obtained by choosing $y(\xi) = 0$.
It is trivial to see that the extremal
$\hat{y}(t) \equiv 0$ we just found
is indeed the global minimizer: $\mathcal{J}(y) \ge 0$
for any function $y$, and $\mathcal{J}(\hat{y}) = 0$.
\end{example}


\subsection{The delta-nabla isoperimetric problem}
\label{subsec:mr:iso}

We consider now delta-nabla isoperimetric problems on time scales.
The problem consists of extremizing
\begin{equation}\label{problem:P:iso}
\mathcal{L}(y) = \left(\int_a^b L_{\Delta}[y](t) \Delta t\right)
\left(\int_a^b L_{\nabla}\{y\}(t) \nabla t\right) \longrightarrow
\textrm{extr}
\end{equation}
in the class of functions $y \in C_{\diamond}^{1}([a,b], \mathbb{R})$
satisfying the boundary conditions
\begin{equation}
\label{bou:con}
y(a) = \alpha \, , \quad y(b) = \beta \, ,
\end{equation}
and the constraint
\begin{equation}
\label{const}
\mathcal{K}(y) = \left(\int_a^b K_{\Delta}[y](t) \Delta t\right)
\left(\int_a^b K_{\nabla}\{y\}(t) \nabla t\right)=k
\end{equation}
where $\alpha$, $\beta$, $k$ are given real numbers.

\begin{definition}
We say that $\hat{y}\in C_{\diamond}^{1}([a,b],
\mathbb{R})$ is a weak local minimizer (respectively weak local
maximizer) for problem \eqref{problem:P:iso}--\eqref{const} if there
exists $\delta
>0$ such that
$$
\mathcal{L}(\hat{y})\leq \mathcal{L}(y) \quad (\text{respectively} \
\   \mathcal{L}(\hat{y}) \geq \mathcal{L}(y))
$$
for all $y \in C_{\diamond}^{1}([a,b], \mathbb{R})$ satisfying the boundary
conditions \eqref{bou:con}, the constraint \eqref{const},
and $||y - \hat{y}||_{1,\infty} < \delta$.
\end{definition}

\begin{definition}
We say that $\hat{y} \in C_{\diamond}^1$ is an extremal for $\mathcal{K}$
if $\hat{y}$ satisfies the delta-nabla integral
equations \eqref{eq:EL1} and \eqref{eq:EL2} for $\mathcal{K}$, \textrm{i.e.},
\begin{multline}
\label{eq:EL1:iso}
\mathcal{K}_\nabla(\hat{y}) \left(\partial_3
K_\Delta[\hat{y}](\rho(t))
-\int_{a}^{\rho(t)} \partial_2 K_\Delta[\hat{y}](\tau) \Delta\tau\right)\\
+ \mathcal{K}_\Delta(\hat{y}) \left(\partial_3
K_\nabla\{\hat{y}\}(t) -\int_{a}^{t} \partial_2
K_\nabla\{\hat{y}\}(\tau) \nabla\tau\right) = \text{const}
\end{multline}
for all $t \in [a,b]_\kappa$;
\begin{multline}
\label{eq:EL2:iso}
\mathcal{K}_\Delta(\hat{y}) \left(\partial_3
K_\nabla\{\hat{y}\}(\sigma(t)) -\int_{a}^{\sigma(t)} \partial_2
K_\nabla\{\hat{y}\}(\tau) \nabla\tau\right)\\
+ \mathcal{K}_\nabla(\hat{y}) \left(\partial_3
K_\Delta[\hat{y}](t)
-\int_{a}^{t} \partial_2 K_\Delta[\hat{y}](\tau) \Delta\tau\right)
= \text{const}
\end{multline}
for all $t \in [a,b]^\kappa$.
An extremizer (\textrm{i.e.}, a weak local minimizer or a weak local
maximizer) for the problem \eqref{problem:P:iso}--\eqref{const} that is
not an extremal for $\mathcal{K}$ is said to be a normal extremizer;
otherwise (\textrm{i.e.}, if it is an extremal for $\mathcal{K}$), the
extremizer is said to be abnormal.
\end{definition}

\begin{theorem}
\label{thm:mr:iso}
If $\hat{y} \in C_{\diamond}^1\left([a,b],\mathbb{R}\right)$
is a normal extremizer for the isoperimetric problem
\eqref{problem:P:iso}--\eqref{const}, then there exists $\lambda \in \mathbb{R}$
such that $\hat{y}$ satisfies the following delta-nabla integral equations:
\begin{multline}
\label{iso:EL1}
\mathcal{L}_\nabla(\hat{y}) \left(\partial_3
L_\Delta[\hat{y}](\rho(t))
-\int_{a}^{\rho(t)} \partial_2 L_\Delta[\hat{y}](\tau) \Delta\tau\right)\\
+ \mathcal{L}_\Delta(\hat{y}) \left(\partial_3
L_\nabla\{\hat{y}\}(t) -\int_{a}^{t} \partial_2
L_\nabla\{\hat{y}\}(\tau) \nabla\tau\right) \\
-\lambda\left\{\mathcal{K}_\nabla(\hat{y}) \left(\partial_3
K_\Delta[\hat{y}](\rho(t))
-\int_{a}^{\rho(t)} \partial_2 K_\Delta[\hat{y}](\tau) \Delta\tau\right)\right.\\
\left.+ \mathcal{K}_\Delta(\hat{y}) \left(\partial_3
K_\nabla\{\hat{y}\}(t) -\int_{a}^{t} \partial_2
K_\nabla\{\hat{y}\}(\tau) \nabla\tau\right)\right\} = \text{const}
\end{multline}
for all $t \in [a,b]_\kappa$;
\begin{multline}
\label{iso:EL2}
\mathcal{L}_\nabla(\hat{y}) \left(\partial_3
L_\Delta[\hat{y}](t)
-\int_{a}^{t} \partial_2 L_\Delta[\hat{y}](\tau) \Delta\tau\right)\\
+ \mathcal{L}_\Delta(\hat{y}) \left(\partial_3
L_\nabla\{\hat{y}\}(\sigma(t)) -\int_{a}^{\sigma(t)} \partial_2
L_\nabla\{\hat{y}\}(\tau) \nabla\tau\right) \\
-\lambda \left\{
\mathcal{K}_\Delta(\hat{y}) \left(\partial_3
K_\nabla\{\hat{y}\}(\sigma(t)) -\int_{a}^{\sigma(t)} \partial_2
K_\nabla\{\hat{y}\}(\tau) \nabla\tau\right)\right.\\
+\left.\mathcal{K}_\nabla(\hat{y}) \left(\partial_3
K_\Delta[\hat{y}](t)
-\int_{a}^{t} \partial_2 K_\Delta[\hat{y}](\tau) \Delta\tau\right)
\right\}= \text{const}
\end{multline}
for all $t \in [a,b]^\kappa$.
\end{theorem}

\begin{proof}
Consider a variation of $\hat{y}$, say $\bar{y}=\hat{y} +
\varepsilon_{1} \eta_{1}+\varepsilon_{2} \eta_{2}$, where for each
$i\in \{1,2\}$, $\eta_{i}\in C_{\diamond}^{1}([a,b],\mathbb{R})$
and $\eta_{i}(a)=\eta_{i}(b)=0$, and $\varepsilon_{i}$ is a
sufficiently small parameter ($\varepsilon_{1}$ and
$\varepsilon_{2}$ must be such that $\parallel
\bar{y}-\hat{y}\parallel_{1,\infty}<\delta$ for some $\delta>0$).
Here, $\eta_{1}$ is an arbitrary fixed function and $\eta_{2}$ is a
fixed function that will be chosen later. Define the real function
\begin{multline*}
\bar{K}(\varepsilon_{1},\varepsilon_{2})
=\mathcal{K}(\bar{y})\\
=\left(\int_a^b
K_{\Delta}[\bar{y}](t) \Delta t\right) \left(\int_a^b
K_{\nabla}\{\bar{y}\}(t) \nabla t\right)-k.
\end{multline*}
We have $\left.\frac{\partial\bar{K}}{\partial
\varepsilon_{2}}\right|_{(0,0)} = 0$, that is,
\begin{multline*}
0 = \mathcal{K}_\nabla(\hat{y})
\int_a^b \left(\partial_2 K_\Delta[\hat{y}](t) \eta_{2}^\sigma(t) +
\partial_3 K_\Delta[\hat{y}](t) \eta_{2}^\Delta(t)\right) \Delta t\\
+ \mathcal{K}_\Delta(\hat{y}) \int_a^b \left(\partial_2
K_\nabla\{\hat{y}\}(t) \eta_{2}^\rho(t) + \partial_3
K_\nabla\{\hat{y}\}(t) \eta_{2}^\nabla(t)\right) \nabla t.
\end{multline*}
Since $\eta_{2}(a)=\eta_{2}(b)=0$,
the first and third integration by parts formula
in \eqref{intBP} give
\begin{multline*}
\int_a^b \partial_2 K_\Delta[\hat{y}](t) \eta_{2}^\sigma(t)\Delta t=
\int_a^t\partial_2 K_\Delta[\hat{y}](\tau)\Delta \tau
\eta_{2}(t)|^{t=b}_{t=a}\\
-\int_a^b\left(\int_a^t
\partial_2 K_\Delta[\hat{y}](\tau)\Delta\tau\right)\eta_{2}^{\Delta}(t) \Delta t\\
=-\int_a^b\left(\int_a^t \partial_2
K_\Delta[\hat{y}](\tau)\Delta\tau\right)\eta_{2}^{\Delta}(t) \Delta t
\end{multline*}
and
\begin{multline*}
\int_a^b \partial_2 K_\nabla\{\hat{y}\}(t) \eta_{2}^\rho(t)\nabla
t=\int_a^t\partial_2 K_\nabla\{\hat{y}\}(\tau)\nabla\tau
\eta_{2}(t)|^{t=b}_{t=a}\\
-\int_a^b\left(\int_a^t\partial_2 K_\nabla\{\hat{y}\}(\tau)\nabla
\tau \right) \eta_{2}^\nabla(t)\nabla t\\
=-\int_a^b\left(\int_a^t\partial_2 K_\nabla\{\hat{y}\}(\tau)\nabla
\tau \right) \eta_{2}^\nabla(t)\nabla t.
\end{multline*}
Therefore,
\begin{multline}
\label{after:parts}
\left.\frac{\partial\bar{K}}{\partial
\varepsilon_{2}}\right|_{(0,0)}
=\mathcal{K}_\nabla(\hat{y}) \int_a^b
\Biggl(\partial_3 K_\Delta[\hat{y}](t)\\
-\int_a^t \partial_2 K_\Delta[\hat{y}](\tau)\Delta\tau\Biggr)
\eta_{2}^{\Delta}(t) \Delta t\\
+\mathcal{K}_\Delta(\hat{y}) \int_a^b \Biggl(\partial_3
K_\nabla\{\hat{y}\}(t)\\
-\int_a^t\partial_2
K_\nabla\{\hat{y}\}(\tau)\nabla \tau \Biggr)
\eta_{2}^\nabla(t)\nabla t.
\end{multline}
Let $$f(t)=\mathcal{K}_\nabla(\hat{y})\left(\partial_3
K_\Delta[\hat{y}](t)-\int_a^t
\partial_2 K_\Delta[\hat{y}](\tau)\Delta\tau\right)$$
and $$g(t)=\mathcal{K}_\Delta(\hat{y})\left(\partial_3
K_\nabla\{\hat{y}\}(t)-\int_a^t\partial_2
K_\nabla\{\hat{y}\}(\tau)\nabla \tau \right).$$ We can then write
equation \eqref{after:parts} in the form
\begin{equation}\label{after:sub}
\left.\frac{\partial\bar{K}}{\partial
\varepsilon_{2}}\right|_{(0,0)}=\int_a^bf(t)\eta_{2}^{\Delta}(t)
\Delta t+\int_a^bg(t)\eta_{2}^\nabla(t)\nabla t.
\end{equation}
Transforming the delta integral in \eqref{after:sub} to a nabla
integral by means of \eqref{eq:DtoN} we obtain
\begin{equation*}
\left.\frac{\partial\bar{K}}{\partial
\varepsilon_{2}}\right|_{(0,0)}=\int_a^bf^{\rho}(t)(\eta_{2}^{\Delta})^{\rho}(t)
\nabla t+\int_a^bg(t)\eta_{2}^\nabla(t)\nabla t
\end{equation*}
and by \eqref{eq:chgN_to_D}
\begin{equation*}
\left.\frac{\partial\bar{K}}{\partial
\varepsilon_{2}}\right|_{(0,0)}=\int_a^b\left(f^{\rho}(t)
+g(t)\right)\eta_{2}^\nabla(t)\nabla t.
\end{equation*}
As $\hat{y}$ is a normal extremizer we conclude, by
Lemma~\ref{DBRL:n} and equation \eqref{eq:EL2:iso}, that there exists
$\eta_2$ such that $\left.\frac{\partial\bar{K}}{\partial
\varepsilon_{2}}\right|_{(0,0)}\neq 0$. Note that the same result
can be obtained by transforming the nabla integral in
\eqref{after:sub} to a delta integral by means of \eqref{eq:NtoD}
and then using Lemma~\ref{DBRL:d} and equation \eqref{eq:EL1:iso}. Since
$\bar{K}(0,0)=0$, by the implicit function theorem we conclude that
there exists a function $\varepsilon_{2}$ defined in the
neighborhood of zero, such that
$\bar{K}(\varepsilon_{1},\varepsilon_{2}(\varepsilon_{1}))=0$, \textrm{i.e.},
we may choose a subset of variations $\bar{y}$ satisfying the
isoperimetric constraint.

Let us now consider the real function
\begin{multline*}
\bar{L}(\varepsilon_{1},\varepsilon_{2})=\mathcal{L}(\bar{y})\\
=\left(\int_a^b
L_{\Delta}[\bar{y}](t) \Delta t\right) \left(\int_a^b
L_{\nabla}\{\bar{y}\}(t) \nabla t\right).
\end{multline*}
By hypothesis, $(0,0)$ is an extremal of $\bar{L}$ subject to the
constraint $\bar{K}=0$ and $\nabla \bar{K}(0,0)\neq \textbf{0}$. By
the Lagrange multiplier rule, there exists some real $\lambda$ such
that $\nabla(\bar{L}(0,0)-\lambda\bar{K}(0,0))=\textbf{0}$. Having
in mind that $\eta_{1}(a)=\eta_{1}(b)=0$, we can write
\begin{multline}
\label{function:L}
\left.\frac{\partial\bar{L}}{\partial
\varepsilon_{1}}\right|_{(0,0)} =\mathcal{L}_\nabla(\hat{y})
\int_a^b \Biggl(\partial_3 L_\Delta[\hat{y}](t)\\
-\int_a^t \partial_2
L_\Delta[\hat{y}](\tau)\Delta\tau\Biggr)\eta_{1}^{\Delta}(t)
\Delta t\\
+\mathcal{L}_\Delta(\hat{y})
\int_a^b \Biggl(\partial_3
L_\nabla\{\hat{y}\}(t)\\
-\int_a^t\partial_2
L_\nabla\{\hat{y}\}(\tau)\nabla \tau \Biggr)
\eta_{1}^\nabla(t)\nabla t
\end{multline}
and
\begin{multline}
\label{function:K}
\left.\frac{\partial\bar{K}}{\partial
\varepsilon_{1}}\right|_{(0,0)}
=\mathcal{K}_\nabla(\hat{y}) \int_a^b
\Biggl(\partial_3 K_\Delta[\hat{y}](t)\\
-\int_a^t \partial_2 K_\Delta[\hat{y}](\tau)\Delta\tau\Biggr)
\eta_{1}^{\Delta}(t) \Delta t\\
+\mathcal{K}_\Delta(\hat{y})
\int_a^b \Biggl(\partial_3
K_\nabla\{\hat{y}\}(t)\\
-\int_a^t\partial_2
K_\nabla\{\hat{y}\}(\tau)\nabla \tau \Biggr)
\eta_{1}^\nabla(t)\nabla t.
\end{multline}
Let $$m(t)=\mathcal{L}_\nabla(\hat{y})\left(\partial_3
L_\Delta[\hat{y}](t)-\int_a^t
\partial_2 L_\Delta[\hat{y}](\tau)\Delta\tau\right)$$
and $$n(t)=\mathcal{L}_\Delta(\hat{y})\left(\partial_3
L_\nabla\{\hat{y}\}(t)-\int_a^t\partial_2
L_\nabla\{\hat{y}\}(\tau)\nabla \tau \right).$$ Then equations
\eqref{function:L} and \eqref{function:K} can be written in the form
\begin{equation*}
\left.\frac{\partial\bar{L}}{\partial
\varepsilon_{1}}\right|_{(0,0)}=\int_a^bm(t)\eta_{1}^{\Delta}(t)
\Delta t+\int_a^bn(t)\eta_{1}^\nabla(t)\nabla t
\end{equation*}
and
\begin{equation*}
\left.\frac{\partial\bar{K}}{\partial
\varepsilon_{1}}\right|_{(0,0)}=\int_a^bf(t)\eta_{1}^{\Delta}(t)
\Delta t+\int_a^bg(t)\eta_{1}^\nabla(t)\nabla t.
\end{equation*}
Transforming the delta integrals in the above equalities to nabla
integrals by means of \eqref{eq:DtoN} and using \eqref{eq:chgN_to_D}
we obtain
\begin{equation*}
\left.\frac{\partial\bar{L}}{\partial
\varepsilon_{1}}\right|_{(0,0)}
=\int_a^b\left(m^{\rho}(t)+n(t)\right)\eta_{1}^\nabla(t)\nabla t
\end{equation*}
and
\begin{equation*}
\left.\frac{\partial\bar{K}}{\partial
\varepsilon_{1}}\right|_{(0,0)}
=\int_a^b\left(f(t)^{\rho}+g(t)\right)\eta_{1}^\nabla(t)\nabla t.
\end{equation*}
Therefore,
\begin{equation}
\label{iso}
\int_{a}^{b}\eta_{1}^{\Delta}(t)\left\{m^{\rho}(t)+n(t)
-\lambda\left(f(t)^{\rho}+g(t)\right)\right\}\nabla t=0.
\end{equation}
Since \eqref{iso} holds for any $\eta_{1}$, by Lemma~\ref{DBRL:n}
we have
\begin{equation*}
m^{\rho}(t)+n(t)-\lambda\left(f(t)^{\rho}+g(t)\right)=c
\end{equation*}
for some $c\in \mathbb{R}$ and all $t \in [a,b]_\kappa$. Hence,
condition \eqref{iso:EL1} holds. Equation \eqref{iso:EL1} can also
be obtained by transforming nabla integrals to delta integrals by
means of \eqref{eq:NtoD} and then using Lemma~\ref{DBRL:d} and
equation \eqref{eq:EL1:iso}.
\end{proof}

In the particular case $L_\nabla \equiv \frac{1}{b-a}$
we get from Theorem~\ref{thm:mr:iso} the main result of
\cite{F:T:09}:

\begin{corollary}[Theorem~3.4 of \cite{F:T:09}]
Suppose that
\begin{equation*}
J(y)=\int_a^b L(t,y^\sigma(t),y^\Delta(t))\Delta t
\end{equation*}
has a local minimum at $y_\ast$
subject to the boundary conditions $y(a)=y_a$ and $y(b)=y_b$
and the isoperimetric constraint
\begin{equation*}
I(y)=\int_a^b g(t,y^\sigma(t),y^\Delta(t))\Delta t = k \, .
\end{equation*}
Assume that $y_\ast$ is not an extremal for the functional $I$.
Then, there exists a Lagrange multiplier constant
$\lambda$ such that $y_\ast$ satisfies the following equation:
\begin{equation*}
\partial_3 F^\Delta(t,y^\sigma_\ast(t),y^\Delta_\ast(t))
-\partial_2 F(t,y^\sigma_\ast(t),y^\Delta_\ast(t))=0
\end{equation*}
for all $t\in[a,b]^{\kappa^2}$,
where $F=L-\lambda g$ and $\partial_3 F^\Delta$
denotes the delta derivative of a composition.
\end{corollary}

One can easily cover abnormal extremizers within our result by
introducing an extra multiplier $\lambda_{0}$.

\begin{theorem}
\label{th:iso:abn}
If $\hat{y} \in C_{\diamond}^1$ is an extremizer for the isoperimetric problem
\eqref{problem:P:iso}--\eqref{const}, then there exist two constants
$\lambda_{0}$ and $\lambda$, not both zero, such that $\hat{y}$
satisfies the following delta-nabla integral equations:
\begin{multline}
\label{iso:EL1:abn}
\lambda_{0}\left\{\mathcal{L}_\nabla(\hat{y})
\left(\partial_3 L_\Delta[\hat{y}](\rho(t))
-\int_{a}^{\rho(t)} \partial_2 L_\Delta[\hat{y}](\tau) \Delta\tau\right)\right.\\
+ \left.\mathcal{L}_\Delta(\hat{y}) \left(\partial_3
L_\nabla\{\hat{y}\}(t) -\int_{a}^{t} \partial_2
L_\nabla\{\hat{y}\}(\tau) \nabla\tau\right) \right\}\\
-\lambda\left\{\mathcal{K}_\nabla(\hat{y}) \left(\partial_3
K_\Delta[\hat{y}](\rho(t))
-\int_{a}^{\rho(t)} \partial_2 K_\Delta[\hat{y}](\tau) \Delta\tau\right)\right.\\
\left.+ \mathcal{K}_\Delta(\hat{y}) \left(\partial_3
K_\nabla\{\hat{y}\}(t) -\int_{a}^{t} \partial_2
K_\nabla\{\hat{y}\}(\tau) \nabla\tau\right)\right\} = \text{const}
\end{multline}
for all $t \in [a,b]_\kappa$;
\begin{multline}
\label{iso:EL2:iso}
\lambda_{0}\left\{\mathcal{L}_\nabla(\hat{y})
\left(\partial_3 L_\Delta[\hat{y}](t)
-\int_{a}^{t} \partial_2 L_\Delta[\hat{y}](\tau) \Delta\tau\right)\right.\\
+\left. \mathcal{L}_\Delta(\hat{y}) \left(\partial_3
L_\nabla\{\hat{y}\}(\sigma(t)) -\int_{a}^{\sigma(t)} \partial_2
L_\nabla\{\hat{y}\}(\tau) \nabla\tau\right)\right\} \\
-\lambda \left\{
\mathcal{K}_\Delta(\hat{y}) \left(\partial_3
K_\nabla\{\hat{y}\}(\sigma(t)) -\int_{a}^{\sigma(t)} \partial_2
K_\nabla\{\hat{y}\}(\tau) \nabla\tau\right)\right.\\
\left. + \mathcal{K}_\nabla(\hat{y}) \left(\partial_3
K_\Delta[\hat{y}](t) -\int_{a}^{t} \partial_2 K_\Delta[\hat{y}](\tau)
\Delta\tau\right) \right\}= \text{const}
\end{multline}
for all $t \in [a,b]^\kappa$.
\end{theorem}

\begin{proof}
Following the proof of Theorem~\ref{thm:mr:iso}, since $(0,0)$ is an
extremal of $\bar{L}$ subject to the constraint $\bar{K}=0$, the
extended Lagrange multiplier rule (see for instance
\cite[Theorem~4.1.3]{Brunt}) asserts the existence of reals
$\lambda_{0}$ and $\lambda$, not both zero, such that
$\nabla(\lambda_{0}\bar{L}(0,0)-\lambda\bar{K}(0,0))=\textbf{0}$.
Therefore,
\begin{equation}\label{iso:abn}
\int_{a}^{b}\eta_{1}^{\Delta}(t)\left\{\lambda_{0}\left(m^{\rho}(t)+n(t)\right)
-\lambda\left(f(t)^{\rho}+g(t)\right)\right\}\nabla t=0.
\end{equation}
Since \eqref{iso:abn} holds for any $\eta_{1}$, by
Lemma~\ref{DBRL:n}, we have
\begin{equation*}
\lambda_{0}\left(m^{\rho}(t)+n(t)\right)-\lambda\left(f(t)^{\rho}+g(t)\right)=c
\end{equation*}
for some $c\in \mathbb{R}$ and all $t \in [a,b]_\kappa$. This
establishes equation \eqref{iso:EL1:abn}. Equation
\eqref{iso:EL2:iso} can be shown using a similar technique.
\end{proof}

\begin{remark}
If $\hat{y} \in C_{\diamond}^1$ is a normal extremizer for the isoperimetric problem
\eqref{problem:P:iso}--\eqref{const}, then we can choose $\lambda_{0}=1$
in Theorem~\ref{th:iso:abn} and obtain Theorem~\ref{thm:mr:iso}. For
abnormal extremizers, Theorem~\ref{th:iso:abn} holds with
$\lambda_{0}=0$. The condition $(\lambda_{0},\lambda)\neq\textbf{0}$
guarantees that Theorem~\ref{th:iso:abn} is a useful necessary
optimality condition.
\end{remark}

In the particular case $L_\Delta \equiv \frac{1}{b-a}$
we get from Theorem~\ref{th:iso:abn} a result of \cite{A:T}:

\begin{corollary}[Theorem~2 of \cite{A:T}]
If $y$ is a local minimizer or maximizer for
\begin{equation*}
I[y]=\int_{a}^{b}f(t,y^\rho(t),y^\nabla(t))\nabla t
\end{equation*}
subject to the boundary conditions $y(a)=\alpha$ and $y(b)=\beta$ and
the nabla-integral constraint
\begin{equation*}
J[y]=\int_{a}^{b}g(t,y^\rho(t),y^\nabla(t))
\nabla t =\Lambda \, ,
\end{equation*}
then there exist two constants $\lambda_0$ and $\lambda$,
not both zero, such that
$$\partial_3 K^\nabla\left(t,y^\rho(t),y^\nabla(t)\right)
-\partial_2 K\left(t,y^\rho(t),y^\nabla(t)\right)=0$$
for all $t \in [a,b]_{\kappa}$,
where $K=\lambda_0 f-\lambda g$.
\end{corollary}

\begin{example}
Let $\mathbb{T}=\{1,2,3,\ldots,M\}$, where $M\in\mathbb{N}$ and
$M\geq 2$. Consider the problem
\begin{equation}\label{ex:iso:1}
\begin{gathered}
\left(\int_{0}^{M}(y^\Delta(t))^2\Delta
t\right)\left(\int_{0}^{M}\left(y^\nabla(t))^2
+y^\nabla(t)\right)\nabla t\right)\rightarrow \min\\
y(0)=0, \quad y(M)=M,
\end{gathered}
\end{equation}
subject to the constraint
\begin{equation}\label{ex:iso:2}
\mathcal{K}(y)=\int_{0}^{M}ty^{\Delta}(t)\Delta t=1.
\end{equation}
Since
\begin{equation*}
\begin{split}
L_{\Delta}&=(y^\Delta)^2,\quad L_{\nabla}=(y^\nabla)^2+y^\nabla,\\
K_{\Delta}&=ty^{\Delta},\quad K_{\nabla}=\frac{1}{M}
\end{split}
\end{equation*}
we have
\begin{equation*}
\partial_2L_{\Delta}=0,\
\partial_3L_{\Delta} =2y^\Delta, \
\partial_2L_{\nabla}=0, \
\partial_3L_{\nabla}=2y^\nabla+1,
\end{equation*}
and
\begin{equation*}
\partial_2K_{\Delta}=0,\quad \partial_3K_{\Delta}=t,
\quad \partial_2K_{\nabla}=0,\quad \partial_3K_{\nabla}=0.
\end{equation*}
As
\begin{multline*}
\mathcal{K}_\Delta(\hat{y}) \left(\partial_3
K_\nabla\{\hat{y}\}(\sigma(t)) -\int_{a}^{\sigma(t)} \partial_2
K_\nabla\{\hat{y}\}(\tau) \nabla\tau\right)\\
+\mathcal{K}_\nabla(\hat{y}) \left(\partial_3K_\Delta[\hat{y}](t)
-\int_{a}^{t} \partial_2 K_\Delta[\hat{y}](\tau) \Delta\tau\right)=t
\end{multline*}
there are no abnormal extremals for the problem
\eqref{ex:iso:1}--\eqref{ex:iso:2}. Applying
equation \eqref{iso:EL2} of Theorem~\ref{thm:mr:iso}
we get the following delta-nabla differential equation:
\begin{equation}\label{ex:iso:3}
2Ay^{\Delta}(t)+B+2By^{\nabla}(\sigma(t))-\lambda t=C,
\end{equation}
where $C\in\mathbb{R}$ and $A$, $B$ are the values of functionals
$\mathcal{L}_{\nabla}$ and $\mathcal{L}_{\Delta}$ in a solution of
\eqref{ex:iso:1}--\eqref{ex:iso:2}, respectively. Since
$y^{\nabla}({\sigma}(t))=y^\Delta(t)$, we can write equation
\eqref{ex:iso:3} in the form
\begin{equation}
\label{ex:iso:4}
2Ay^{\Delta}(t)+B+2By^\Delta-\lambda t=C.
\end{equation}
Observe that $B\neq 0$ and $A>2$. Hence, solving equation
\eqref{ex:iso:4} subject to the boundary conditions $y(0)=0$ and
$y(M)=M$ we get
\begin{equation}\label{ex:iso:5}
y(t)= \left[1 - \frac{\lambda\left(M-t\right)}{4(A+B)}\right] t \, .
\end{equation}
Substituting \eqref{ex:iso:5} into \eqref{ex:iso:2} we obtain
$\lambda=-\frac{\left( A+B \right) \left( M-2 \right)}{12 M
\left(M-1 \right)}$. Hence,
\begin{equation*}
y(t)=\frac{\left(4\,{M}^{2}-7\,M-3 M \,t + 6\,t\right) t}{M
\left( M-1\right)}
\end{equation*}
is a candidate local minimizer for the problem
\eqref{ex:iso:1}--\eqref{ex:iso:2}.
\end{example}


\section{Conclusion}

The calculus of variations on time scales is an important subject
under strong current research (see
\cite{MyID:153,MyID:140,MyID:159,MyID:141,MyID:183,Basia:post_doc_Aveiro:3}
and references therein). Here we review
a general necessary optimality condition
for problems of the calculus of variations
on time scales \cite{delfim:Bedlewo:2009,MOMA09}.
The proposed calculus of variations
extends the problems with delta derivatives
considered in \cite{B:04,F:T:09} and analogous nabla
problems \cite{A:T,NM:T} to more general cases
described by the product of a delta and a nabla
integral. Minimization of functionals
given by the product of two integrals
were considered by Euler himself,
and are now receiving an increasing interest
because of their nonlocal properties
and their applications in economics \cite{Pedregal}.
We proved Euler--Lagrange type conditions
for the generalized calculus of variations
and corresponding natural boundary conditions.

The calculus of variations here promoted
can be further developed.
For instance, we can continue our study
by proving sufficient optimality conditions.
Moreover, the results here presented can be generalized in different ways:
(i) to variational problems
involving higher-order delta and nabla derivatives,
unifying and extending the higher-order results on time scales
of \cite{F:T:08} and \cite{NM:T};
(ii) to problems of the calculus of variations on time scales
introduced in \cite{MyID:151}, with a functional
which is the composition of a certain scalar
function $H$ with the delta integral of a vector valued field $f_\Delta$
and a nabla integral of a vector field $f_\nabla$, \textrm{i.e.}, of the form
$$
H\left(\int_{a}^{b}f_\Delta(t,y^{\sigma}(t),y^{\Delta}(t))\Delta t,
\int_{a}^{b}f_\nabla(t,y^{\rho}(t),y^{\nabla}(t))\nabla t\right).
$$
Euler--Lagrange equations and natural boundary conditions for such
problems on time scales can be proved,
and present results obtained as corollaries.


\subsection*{Acknowledgements}

This work is supported by the \emph{Portuguese Foundation
for Science and Technology} (FCT) through the
\emph{Center for Research and Development
in Mathematics and Applications} (CIDMA).



\end{multicols}{}

\end{document}